\documentclass[12pt]{amsart}
\usepackage{graphicx}
\usepackage[linktocpage=true,colorlinks,citecolor=blue,linkcolor=blue,urlcolor=blue]{hyperref}
\usepackage{amssymb}
\usepackage{cite}
\usepackage{amsmath}
\usepackage{latexsym}
\usepackage{amscd}
\usepackage{tikz}
\usepackage{amsthm}
\usepackage{mathrsfs}
\usepackage{url}
\usepackage[utf8]{inputenc}
\usepackage[english]{babel}
\usepackage{amsfonts}
\usepackage{mathtools}

\def\P{\mathbb P}

\def\E{\mathbb E}

\newtheorem{theorem}{Theorem}[section]
\newtheorem{lemma}[theorem]{Lemma}

\theoremstyle{remark}
\newtheorem{remark}[theorem]{Remark}

\theoremstyle{definition}

\theoremstyle{remark}

\numberwithin{equation}{section}

\begin{document}
	\title{ Oscillations of random multiplicative functions under initial bias}
\author{ Rodrigo Angelo }
\address{São Paulo, SP, Brazil}
\email{rsa365@gmail.com}
\author{Max Wenqiang Xu}
\address{Yau Mathematical Sciences Center, Tsinghua University, Beijing, China}
\email{maxxu1729@gmail.com}
	\begin{abstract}
We prove that if $f$ is a random completely multiplicative function, conditional on $f(p)=1$ for each prime $p \le (\log x)^{2-\epsilon}$, the probability that $\sum_{1\le n \le N}f(n)\ge 0$ for all $N\le x$ is $o(1)$ as $x \rightarrow \infty$. This solves a conjecture of Kucheriaviy, who has a complementary result showing this exponent is sharp.

We also prove that almost surely the partial sums of $\sum\frac{f(n)}{\sqrt{n}}$ change signs infinitely many times, solving a problem of Aymone.
\end{abstract}
	\maketitle
\section{Introduction}
Let $f$ denote the random completely multiplicative function, defined by $f(p) = \pm 1$ with probability $\frac{1}{2}$ independently at each prime $p$, and extended completely multiplicatively to all positive integers $n$. In this paper we solve two problems on sign changes of partial sums of this function.

The first problem concerns the effect of the initial primes on these sign changes. Kucheriaviy \cite{kucheriaviy2025} proved that, conditional on $f(p) = 1$ for each $p \le y$, the probability that $\sum_{n \le N}f(n) \ge 0$ for each $N \le x$ is $1-o(1)$, as long as $y \ge C (\log x)^2 \frac{\log \log x}{\log \log \log x}$ for a constant $C$. One explanation for the scale $y= (\log x)^{2+o(1)}$ is that at this threshold the number of $y$-smooth integers up to $x$ surpasses $x^{\frac{1}{2}+o(1)}$, which is the square root of the variance of the remaining random sum. We prove a complementary result confirming his Conjecture~1 in \cite{kucheriaviy2025}.

\begin{theorem}\label{maintheorem}
Let $f$ be the random completely multiplicative function. Let $x$ be large and suppose $y = o\left(\left(\frac{\log x}{\log \log x}\right)^2\right)$. Then, conditional on $f(p)=1$ for each $p \le y$, the probability that $\sum_{n \le N} f(n) \ge 0$ for each $N \le x$ is $o(1)$ as $x \rightarrow \infty$.
\end{theorem}

The second problem was proposed by Aymone \cite{aymone2024ii}, on sign changes of the partial sums of $\sum\frac{f(n)}{n^{\sigma}}$. For $\sigma > \frac{1}{2}$ this sum almost surely converges to its Euler product, which is a positive number, so the partial sums change signs only finitely many times. For $\sigma < \frac{1}{2}$, he established that almost surely this sum changes sign infinitely often. We solve the remaining case.\footnote{This problem had been solved by the first author in an earlier preprint that was incorporated into this paper.}

\begin{theorem}\label{thm: sign changes}
    Let $f$ be a sample of the random completely multiplicative function. Almost surely the partial sums of $\frac{f(n)}{\sqrt{n}}$ change signs infinitely many times. 
\end{theorem}

A deterministic variation is studied in \cite{betweenpolyaandturan}, where it is suggested that the partial sums of $\sum \frac{\lambda(n)}{\sqrt{n}}$ may be eventually always negative, where $\lambda$ is the Liouville function.  Our proof relies on oscillations near $\frac{1}{2}^{+}$ specific to random multiplicative functions, so we don't think it affects the plausibility of finitely many sign changes in this deterministic version, which is certainly very subtle since the eventual negativity of this sum would imply both the Riemann Hypothesis and the simplicity of the zeros of $\zeta(s)$. See also Remark \ref{notpositive}, where we show the partial sums of $\sum \frac{\lambda(n)}{\sqrt{n}}$ may not be eventually always positive.

These problems are relevant to a question of much interest to us --- the probability that $f$ has nonnegative partial sums up to $x$.  A consequence of Kucheriaviy's result on $f$ with initial bias is that this probability is at least $\exp(-C\frac{(\log x)^2}{\log \log \log x})$.  Determining the asymptotic decay of this probability is an open problem, and currently this is the best known lower bound. Kalmynin \cite{kal2021} gives an upper bound of the shape $(\log x)^{-c}$ for this probability. A deterministic analogue of this question is the enumeration of quadratic characters $\chi$ with the property that all of their partial sums are positive. It seems likely that there are infinitely many such characters, but that is an open question. In an upcoming collaboration with Soundararajan \cite{ASX} we prove an upper bound of the shape $(\log x)^{-c}$ for the probability that a quadratic character has this property, among characters of discriminant at most $x$.

These upper bounds are based on the fact that the nonnegativity constrains the maximum absolute value of $\sum \frac{f(n)}{n^s}$ over complex $s$ in a fixed vertical interval near the critical line. Aymone also used this approach for establishing the case $\sigma < \frac{1}{2}$ of Theorem \ref{thm: sign changes}. We note Aymone and Kucheriaviy used a result of Harper on the suprema of Gaussian processes \cite{Harpersup} as a key input. A similar maximum on the critical line is also the object of study of the Fyodorov--Hiary--Keating conjecture \cite{FHK}, which has been the subject of a lot of recent work, and is now a theorem \cite{fhk1, fhk2}.

While this approach with complex $s$ has proven effective in some cases, it doesn't seem to work in the problems we present here. Instead, we use consequences of nonnegativity only to real valued $s$, namely that it implies roughly decreasing behavior of the series $\sum \frac{f(n)}{n^s}$, and we develop a minimal approach for handling the joint behavior of this series at various $s\rightarrow \frac{1}{2}^+$ to show this decreasing behavior is atypical, an argument more akin to \cite{BM89} or  \cite{klurman2024signchangesshortcharacter}.

See also \cite{lamzouri2026realzeroslschid} for additional results on oscillations of quadratic characters, and \cite{Sarnak2011Letter, Jung2016} for adjacent work in the context of automorphic forms.

\subsection*{Acknowledgement} 
M.W.X. was supported by a Simons Junior Fellowship from the Simons Foundation.

\section{Proof of Theorem~\ref{maintheorem}}

Throughout this paper, we set $f(p) = 1$ for primes $p \le y$, while $f(p) = \pm 1$ with probability $\frac{1}{2}$ for $p>y$. In the second problem, where there is no conditioning on $p \le y$, the same lemmas may be reused by simply setting $y=1$. A summation over $p$ always denotes a sum over primes, while a summation over $n$ is always over the positive integers.

In both problems, we relate the nonnegativity of the sums to approximately decreasing behavior of the series $\sum_n \frac{f(n)}{n^{\frac{1}{2}+t}}$ for small $t>0$. With the Euler product this becomes the exponential of a sum of independent random variables, namely $\sum \frac{f(n)}{n^{\frac{1}{2}+t}} = \exp(R(t)+O(1))$ where

\begin{equation}\label{eqn: R}
R(t) = \sum_{p} \frac{f(p)}{p^{\frac{1}{2}+t}} + \frac{1}{2}\sum_{p} \frac{1}{p^{1+2t}}.   
\end{equation}

A central limit theorem for the differences $R(t)-R(2t)$ will generate oscillations that are ultimately inconsistent with the decreasing behavior. Notice that $R(t)$ is like the ``$\log$" of $\sum_n \frac{f(n)}{n^{\frac{1}{2}+t}}$, while $R(t)-R(2t)$ is like its derivative, so these differences serve as an alternative to looking at a log derivative of this series.

\begin{lemma} \label{centrallimit}
Let $K$ be a fixed positive integer. Let $y \ge 1$ and $0 < t_1<\dots<t_K$. Let $f$ be the random completely multiplicative function conditioned on $f(p) = 1$ for $p \le y$. Then as long as $\frac{t_{i+1}}{t_i}\rightarrow \infty$ for each $i$, and $y^{\frac{1}{2}}t_K \rightarrow 0$, we have the following limit in distribution of random variables
$$
\left( R(t_1)-R(2t_1), R(t_2)-R(2t_2), \dots, R(t_{K})- R(2t_{K}) \right) \xrightarrow{d} (N_1, N_2, \dots, N_K),
$$
where $(N_1,\dots,N_K)$ is a joint normal distribution with independent components $N_i$, each with mean $\frac{\log 2}{2}$ and variance $\log \frac{9}{8}$.
\end{lemma}

\begin{proof}[Proof of Lemma \ref{centrallimit}]
We have
    $$R(t)-R(2t) = \sum_{p \le y} \left(\frac{1}{p^{\frac{1}{2}+t}}-\frac{1}{p^{\frac{1}{2}+2t}}\right) + \sum_{p>y} f(p)\left(\frac{1}{p^{\frac{1}{2}+t}}-\frac{1}{p^{\frac{1}{2}+2t}}\right) +\frac{1}{2}\sum_{p}\left(\frac{1}{p^{1+2t}}-\frac{1}{p^{1+4t}}\right).$$
    The first sum, coming from $p \le y$, is bounded by 
    $$
    \ll t \sum_{p \le y} \frac{\log p}{p^{\frac{1}{2}}} \ll t y^{\frac{1}{2}}.$$
    For each $t = t_i$, this is $o(1)$ since $y^{\frac{1}{2}} t_i \le y^{\frac{1}{2}} t_K= o(1)$.
    The contribution of the third sum is, via Lemma \ref{approximation}, 
    $\frac{1}{2} (\log\frac{1}{2t} - \log\frac{1}{4t}) + o(1) = \frac{\log 2}{2} + o(1)$. This constant gives the means of our Gaussian. 
    
    To achieve our Lemma we only need to compute the covariance matrix of the limiting distribution. For the diagonal terms, an application of Lemma~\ref{approximation} gives $$\sigma(t_i,t_i) = \sum_{p>y} \left(\frac{1}{p^{\frac{1}{2}+t_i}}-\frac{1}{p^{\frac{1}{2}+2t_i}}\right)^2 = \log\frac{1}{4t_i} + \log\frac{1}{2t_i}-2 \log \frac{1}{3t_i}+o(1) = \log\frac{9}{8} + o(1).$$
    And again Lemma \ref{approximation} implies that the covariances are
    \begin{equation}\label{eqn: covariance}
       \sigma(t_i,t_j) = \sum_{p>y} \left(\frac{1}{p^{\frac{1}{2}+t_i}}-\frac{1}{p^{\frac{1}{2}+2t_i}}\right)\left(\frac{1}{p^{\frac{1}{2}+t_j}}-\frac{1}{p^{\frac{1}{2}+2t_j}}\right) = \log \frac{t_i+2t_j}{t_i+t_j}+\log\frac{2t_i+t_j}{2t_i+2t_j}+o(1)=o(1). 
    \end{equation}
   Here we have used that for $i>j$, $\frac{t_i}{t_j} \rightarrow \infty$.
   
    Since $\sup_{p, t_i}\left| \frac{1}{p^{\frac{1}{2}+t_i}}-\frac{1}{p^{\frac{1}{2}+2t_i}}\right|\rightarrow 0$, the completion of our proof is just an application of a $K$-dimensional central limit theorem to the term in the middle. For example we may use the version of a Berry--Esseen bound from \cite{lyapunovbound}. They show that if $X_n$ is a sequence of independent $K$-dimensional random variables with $\E[X_n] = 0$, and $X = \sum_n X_n$ then 
    
    $$\sup_{R} |\P[X \in R] - \P[N \in R]| \le cK^{\frac{1}{4}} \sum_{n} \E|X_n|^3$$
    for an absolute constant $c$, where the supremum is over all rectangles $R$ in $\mathbb{R}^K$, and $N$ is the normal distribution with covariance matrix equal to the covariance matrix of $X$. While their result is only stated for finite sums, it is clear that as long as the partial sums of $X_n$ converge in distribution to $X$, the bound extends to infinite sums. In our case this is guaranteed by the fact that $\sum \E|X_n|^2$ converges.

    We use the random vectors $X_p = \left(f(p)\left(\frac{1}{p^{\frac{1}{2}+t_{i}}}-\frac{1}{p^{\frac{1}{2}+2t_{i}}}\right)\right)_{i=1}^K$. Their bound becomes

    $$cK^{\frac{1}{4}} \sum_{p>y} \E|X_p|^3 \le c K^{\frac{1}{4}}\sum_{p>y}\sup_{t_i}\left| \frac{1}{p^{\frac{1}{2}+t_i}}-\frac{1}{p^{\frac{1}{2}+2t_i}}\right|\E|X_p|^2.$$

    Our variance computation shows $\sum_{p>y} \E|X_p|^2$ converges to a quantity that tends to $K\log \frac{9}{8}$, and evidently $\sup_{t_i}\left| \frac{1}{p^{\frac{1}{2}+t_i}}-\frac{1}{p^{\frac{1}{2}+2t_i}}\right|\rightarrow 0$ uniformly on $p$, since $t_K \rightarrow 0$. Therefore this sup distance between $\sum X_p$ and the Gaussian random variable $N'$ with covariance matrix $\sigma(t_i, t_j)$ tends to zero --- this is equivalent to convergence in distribution to zero of $\sum X_p - N'$. Our computations show this covariance matrix converges to a diagonal matrix of values $\log \frac{9}{8}$, which shows this random variable $N'$ converges in distribution to $K$ independent copies of Gaussians of variance $\log \frac{9}{8}$, which completes the proof.

\end{proof}

The following was the only estimate needed. Keep in mind that $t \log 2y \ll y^{\frac{1}{2}}t \rightarrow 0$, so it applies for $t = t_i$ under the assumptions of Lemma \ref{centrallimit}.

\begin{lemma}\label{approximation} Let $y\ge 1$ and $t >0$. There is a constant $c$ such that if $t \log 2y  \rightarrow 0$, then 
    $$\sum_{p>y} \frac{1}{p^{1+t}} = \log \frac{1}{t}+c-\sum_{p \le y} \frac{1}{p}+ o(1).$$
\end{lemma}

\begin{proof}
    Simply write $\sum_{p>y} \frac{1}{p^{1+t}} = \sum_{p} \frac{1}{p^{1+t}} - \sum_{p\le y} \frac{1}{p^{1+t}}$. The estimate $\sum_{p} \frac{1}{p^{1+t}} = \log \frac{1}{t} +c+o(1)$ is classical from looking at $\log \zeta(1+t)$, while
    $$\sum_{p\le y} \left(\frac{1}{p}- \frac{1}{p^{1+t}}\right) \ll \sum_{p \le y} t\frac{\log p}{p} \ll t \log 2y = o(1).$$
\end{proof}

Now we look at how nonnegative partial sums constrain
$$F(s) = \sum_{n \le x} \frac{f(n)}{n^s}.$$

\begin{lemma}\label{decreasing} 
    If $f$ has nonnegative partial sums up to $x$, then $\frac{F(s)}{s}$ is decreasing for $s>0$.
\end{lemma}

\begin{proof}
    With integration by parts we may write
    $$F(s) = \sum_{n \le x} \frac{f(n)}{n^s} = \int_{1}^x \sum_{n \le z} f(n) \frac{s}{z^{s+1}}dz +\frac{\sum_{n\le x} f(n)}{x^s}.$$

    Divide both sides by $s$ and differentiate to obtain
    $$\frac{d}{ds}\frac{F(s)}{s} = \int_{1}^x \sum_{n \le z} f(n) \frac{d}{ds}\frac{1}{z^{s+1}}dz +\sum_{n\le x} f(n)\frac{d}{ds}\frac{1}{sx^s}.$$
    
Because $\frac{1}{z^{s+1}}$ over $1 \le z \le x$ and $\frac{1}{sx^s}$ are both decreasing functions of $s>0$, and the partial sums of $f$ are nonnegative, this implies $\frac{d}{ds}\frac{F(s)}{s}\le 0$ for $s>0$, as desired.
\end{proof}

For $s = \frac{1}{2}+t$ appropriately above $\frac{1}{2}$, $F(s)$ may be relaxed to the infinite sum, which has the Euler product.  With this Euler product we write $F(s)$ as the exponential of the sum over primes $R(t)$, which recall we have a central limit theorem for.

\begin{lemma}\label{transition} 
     Let $f$ denote the random completely multiplicative function, conditioned on $f(p) = 1$ for $p \le y$. Assume $\frac{20\log \log x}{\log x}\le t \le \frac{1}{y^{\frac{1}{2}}}$. Then for $s = \frac{1}{2} + t$, with probability $1-o(1)$ one has
$$\sum_{n \le x} \frac{f(n)}{n^s} = \exp(R(t)+O(1)).$$
\end{lemma}

\begin{proof}
We begin by showing
$$
\mathbb{E} \left| \sum_{n \leq x} \frac{f(n)}{n^s}  \prod_p \left( 1 - \frac{f(p)}{p^s} \right) -1\right|^2 =o(1).
$$
By writing the finite sum as the entire sum minus its tail, and using the Euler product we get that this is equal to
$$\mathbb{E} \left| \sum_{n > x} \frac{f(n)}{n^s}  \prod_{p} \left( 1 - \frac{f(p)}{p^s} \right) \right|^2.$$ 
    Factoring out $p \le y$ gives that
\begin{equation}\label{eqn: full}
    \prod_{p \le y} \left(1  - \frac{1}{p^s}\right)^2\mathbb{E} \left| \sum_{n > x} \frac{f(n)}{n^s}  \prod_{p>y} \left( 1 - \frac{f(p)}{p^s} \right) \right|^2.    
    \end{equation}
Now notice that $\E[f(mn)]=1$ precisely when $mn$ is of the shape $ab^2$ where all the prime factors of $a$ are $\le y$ and the prime factors of $b$ are $>y$. Otherwise $\E[f(mn)] = 0$.
    
    Expand the inner sum and consider the coefficient of $f(N)$ for each $N = ab^2$ as described. We can write $ab^2=n_1n_2m_1m_2$, where $n_1,n_2$ are from the sum $\sum_{n>x}\frac{f(n)}{n^{s}}$ and $m_1,m_2$ from the infinite product. Since the prime factors of $a$ are all $\le y$, we may further write $n_1=d_1\ell_1$ and $n_2=d_2\ell_2$ where $d_1d_2=a$, and $m_1m_2\ell_1\ell_2=b^2$. As such, the number of times $f(N)$ with $N=ab^2$ appears in the expansion is at most $\tau_2(a)\tau_4(b^2)$.
Together with an application of Rankin's trick and $s = \frac{1}{2}+t$, this implies that  $\mathbb{E} \left| \sum_{n > x} \frac{f(n)}{n^s}  \prod_{p>y} \left( 1 - \frac{f(p)}{p^s} \right) \right|^2$
\begin{equation*}
    \begin{split}
    \le \sum_{\substack{ab^2>x^2 \\ p|a \implies p\le y\\ p|b \implies p>y}} \frac{\tau_2(a)\tau_4(b^2)}{(ab^2)^s} &\ll \frac{1}{x^t}\sum_{\substack{p|a \implies p\le y\\ p|b \implies p>y}} \frac{\tau_2(a)\tau_4(b^2)}{(ab^2)^{\frac{1}{2}+\frac{t}{2}}} \\ & \ll \frac{1}{x^t}\prod_{p \le y} \left(1-\frac{1}{p^{\frac{1}{2}+\frac{t}{2}}}\right)^{-2} \prod_{p>y} \left(1-\frac{1}{p^{1+t}}\right)^{-10}.     
    \end{split}
\end{equation*}
    Given this, the contribution of all primes $p \le y$ in \eqref{eqn: full} is 
    $$\prod_{p \le y} \left(1  - \frac{1}{p^{\frac{1}{2}+t}}\right)^2\left(1-\frac{1}{p^{\frac{1}{2}+\frac{t}{2}}}\right)^{-2} \ll \exp \left(t\sum_{p \le y} \frac{\log p}{p^{\frac{1}{2}}}\right) \ll \exp(t y^{\frac{1}{2}}).$$
 The contribution from primes $p>y$ is 
    $$\prod_{p>y} \left(1-\frac{1}{p^{1+t}}\right)^{-10} \le \zeta(1+t)^{10} \ll \frac{1}{t ^{10}}.$$
    Therefore, \eqref{eqn: full} is 
    $$\ll \frac{\exp(ty^{\frac{1}{2}})}{x^{t}t^{10}}=o(1).$$
    Note that $x^{t} \gg (\log x)^{20}$, $t^{10} \ge \frac{1}{(\log x)^{10}}$, and $\exp(ty^{\frac{1}{2}}) \ll 1$,  using $\frac{20\log \log x}{\log x}\le t \le \frac{1}{y^{\frac{1}{2}}}$.

  Given the above second moment estimate, we can see that with probability at least $1 - o(1)$ one has
    $$\frac{1}{2}\prod_{p} \left( 1 - \frac{f(p)}{p^s} \right)^{-1} \le \sum_{n \le x} \frac{f(n)}{n^s} \le 2\prod_{p} \left( 1 - \frac{f(p)}{p^s} \right)^{-1}.$$
    The lemma then follows from the fact that
    $$\prod_{p} \left( 1 - \frac{f(p)}{p^s} \right)^{-1} = \exp(R(t)+O(1)),$$
via the approximation $\left(1-\frac{f(p)}{p^s}\right)^{-1} = \exp\left(\frac{f(p)}{p^s}+\frac{1}{2p^{2s}}+O\left(\frac{1}{p^\frac{3}{2}}\right)\right)$.

\end{proof}

\iffalse

\begin{lemma}\label{decreasing} 
    Let $t>0$. If $f$ has nonnegative partial sums up to $x$, then $F(\frac{1}{2}+2t) \le F(\frac{1}{2}+t)+O(t)$.
\end{lemma}

\begin{proof}

Let $f(n)$ have nonnegative partial sums up to $x$. Write
    $$F(s) = \sum_{n \le x} \frac{f(n)}{n^s} = \sum_{n \le x} (f(1)+\dots+f(n))\left(\frac{1}{n^s}-\frac{1}{(n+1)^s}\right) + \frac{f(1)+\dots+f(x)}{(x+1)^s}.$$
    It is not hard to verify that functions $W_n(s) = \frac{1}{n^s}- \frac{1}{(n+1)^s}$ are decreasing functions of $s \ge \frac{1}{2}$ as long as $n \ge 7$ (indeed $W_n'(\frac{1}{2}) = \frac{\log(n+1)}{\sqrt{n+1}} - \frac{\log n}{\sqrt{n}}$, which is negative for $n \ge 7$). The function $V_n(t) = \frac{1}{(n+1)^t}$ is clearly decreasing as well.

    Only the contribution of $n \le 6$ is not decreasing, but clearly
    
    $$\sum_{n \le 6} (f(1)+\dots+f(n))\left(\frac{1}{n^{\frac{1}{2}+t}}-\frac{1}{(n+1)^{\frac{1}{2}+t}}\right)  = O(t),$$
which is enough to complete the proof.

\end{proof}

\fi

\begin{proof}[Proof of Theorem \ref{maintheorem}]

Fix $K>0$ and choose
$$t_i = \frac{\log \log x}{ \log x} A^i$$ for $i = 1,\dots, K$. If $y = o \left(\frac{\log x}{\log \log x}\right)^2$, we may choose $A \rightarrow \infty$ growing slowly enough so that
$$yA^{2K} = o\left(\frac{\log x}{\log \log x}\right)^2$$ as well, so Lemma ~\ref{centrallimit} is applicable for $R(2t_i)-R(t_i)$ as $x \rightarrow \infty$.

Suppose that $f$ has nonnegative partial sums $\sum_{1\le n \le N}f(n)$ for $ N\le x$. Lemma~\ref{decreasing} implies that for $t>0$
$$\frac{F(\frac{1}{2}+2t)}{2} \le F(\frac{1}{2}+t).$$
After applying Lemma \ref{transition} we obtain an absolute constant $C$ such that $$R(2t_{i})-R(t_i) \le C$$ for all $1\le i \le K$ with exceptional probability $o_{x\to \infty}(1)$.

By using Lemma~\ref{centrallimit} we bound the probability that $R(2t_{i})-R(t_i) \le C$ for all $1\le i \le K$ by the probability that $N_i \ge -C$ for each such $i$, plus an error that is $o_{x\to \infty}(1)$. From independence, this probability is simply $$(1-\delta)^K$$ for some positive $\delta>0$, where $1-\delta$ denotes the probability that $N_i \ge -C$. By making $K$ arbitrarily large, we conclude that the probability that $f$ has nonnegative partial sums up to $x$ must tend to $0$ as $x \rightarrow \infty$, as desired.

\end{proof}
\begin{remark}  
A note about the final scale of $t$ and $y$: for $t$ beyond $\frac{1}{y^{\frac{1}{2}}}$, the contribution of $p \le y$ to the value of $R(t)-R(2t)$ becomes too substantial, and their behavior is typically solidly decreasing regardless of the values of $f(p)$ for $p>y$ (we also lose the ability of relating $R(t)$ to the Euler product in this range). Additionally, the transition from $\sum_{n \le x} \frac{f(n)}{n^{\frac{1}{2}+t}}$ to the Euler product can only be achieved for $t \gg \frac{\log \log x}{\log x}$ --- this is what limits how large $y$ can be.

Notice also that most of the contribution to $R(t)-R(2t)$ comes from primes around $\exp(\frac{1}{t})$, and our final choice of $t$ is  $\gg \frac{\log \log x}{\log x}$, so the variance coming from primes below $\exp(o(\frac{\log x}{\log \log x}))$ is never used! Only the randomness of primes above it is needed to the oscillations we generate.

\end{remark}

\section{Proof of Theorem \ref{thm: sign changes}}

\begin{lemma} \label{almostdecreasing}
Let $f(n)$ be a random completely multiplicative function. For $t>0$, define 
 $$S(t) = \frac{1}{t} \sum_{n=1}^{\infty} \frac{f(n)}{n^{\frac{1}{2}+t}}.$$
If $\sum_{n \le N} \frac{f(n)}{\sqrt{n}} \ge 0$ holds for all $N \ge N_0$, then for any $0<t_1<t_2<\frac{1}{4}$ one has
$$S(t_2)\le S(t_1)+8 \sqrt{N_0}.$$
Similarly, if $\sum_{n \le N} \frac{f(n)}{\sqrt{n}} \le 0$ holds for all $N \ge N_0$, then for any $0<t_1<t_2<\frac{1}{4}$ one has
$$S(t_1)\le S(t_2)+8 \sqrt{N_0}.$$
\end{lemma}

\begin{proof}
Integrating by parts
    $$S(t) = \frac{1}{t}\sum_{n=1}^{\infty}\frac{f(n)}{n^{\frac{1}{2}+t}} = \int_{1}^\infty \sum_{n \le z} \frac{f(n)}{\sqrt{n}} \frac{1}{z^{t+1}}dz.$$
Since $\frac{1}{z^{t+1}}$ is decreasing, we obtain that for $t_1<t_2$, it holds that
\begin{equation*}
    \begin{split}
 S(t_2)-S(t_1) &\le \int_{1}^{N_0} \sum_{n \le z} \frac{f(n)}{\sqrt{n}} \left(\frac{1}{z^{t_2+1}}-\frac{1}{z^{t_1+1}}\right)dz\\& \le  \int_1^{N_0} 2\sqrt{z}\frac{1}{z^{t_1+1}} \le 2\frac{N_0^{\frac{1}{2}-t_1}}{\frac{1}{2}-t_1} \le  8 \sqrt{N_0}  .     
    \end{split}
\end{equation*}
The other case can be proved the same way.
\end{proof}

We are now ready to prove Theorem~\ref{thm: sign changes}.
\begin{proof}[Proof of Theorem \ref{thm: sign changes}]
For a given value of $N_0$, we will show that the event that $\sum_{n \le N} \frac{f(n)}{\sqrt{n}} \ge 0$ for all $N \ge N_0$ has probability $0$.

By using the almost sure Euler product, we have the approximation
$$S(t) = \exp\Big(\log \frac{1}{t}+R(t)+O(1)\Big)$$
where as previously defined $R(t) = \sum_{p} \frac{f(p)}{p^{\frac{1}{2}+t}} + \frac{1}{2}\sum_{p} \frac{1}{p^{1+2t}}$.

By applying Lemma \ref{almostdecreasing}, we see that if $\sum_{n \le N} \frac{f(n)}{\sqrt{n}} \ge 0$ for each $N \ge N_0$ we have $S(2t) \le S(t)+8\sqrt{N_0}$. 
The event $S(t)<8\sqrt{N_0}$ has probability $o(1)$ as $t \rightarrow 0$.  This is because $\sum_{p} \frac{f(p)}{p^{\frac{1}{2}+t}} \le - \log \frac{1}{t}$ holds with probability $\ll\frac{1}{\log \frac{1}{t}}$, which can be proved by using $\E[|\sum_{p} \frac{f(p)}{p^{\frac{1}{2}+t}}|^{2}] \le \log \frac{1}{t} +O(1)$.

Given the event $S(t)\ge8\sqrt{N_0}$ holds, the inequality $S(2t) \le S(t)+8\sqrt{N_0}$ implies that $S(2t) \le 2S(t)$, which in turn implies $R(2t) \le R(t)+C$ for some constant $C$. Apply Lemma \ref{centrallimit}, setting $y = 1$ and say $t_1,\dots, t_K = t^K, t^{K-1},\dots, t$ for some $1>t>0$. We obtain values $t_1,\dots, t_K$ such that the distribution of $(R(t_i)-R(2t_i))_{1\le i \le K}$ are close to independent Gaussians $N_i$. By making $t$ tend to $0$ we obtain that the probability $R(2t) -R(t) \le C$ is at most the probability that $N_i \ge -C$ for each $i = 1,\dots, K$, which is just $(1-\delta)^K$ for an absolute $\delta >0$.

Because $K$ can be made arbitrarily large, this means that the probability $\sum_{n \le N} \frac{f(n)}{\sqrt{n}} \ge 0$ for every $N \ge N_0$ is $0$. Analogously, the probability $\sum_{n \le N} \frac{f(n)}{\sqrt{n}} \le 0$ for each $N \ge N_0$ is 0. The event that $\sum_{n\le N} \frac{f(n)}{\sqrt{n}}$ changes signs finitely many times is the union of these countably many events indexed by $N_0 \ge 1$, and therefore has probability 0 as well.

\end{proof}

\begin{remark} \label{notpositive}
    For the deterministic sums $\sum_{n \le x} \frac{\lambda(n)}{\sqrt{n}}$, the argument from Lemma \ref{almostdecreasing} implies that this sum can't be positive for every sufficiently large $x$.
    
    It is shown in \cite{betweenpolyaandturan} that if these partial sums change signs finitely many times, the Riemann Hypothesis follows. Under the Riemann Hypothesis, $S(t) = \frac{1}{t} \sum \frac{\lambda(n)}{n^{\frac{1}{2}+t}}$ converges to $\frac{1}{t} \frac{\zeta(1+2t)}{\zeta(\frac{1}{2}+t)}\sim \frac{1}{2t^2 \zeta(\frac{1}{2})}$. But if the partial sums were positive for $x \ge N_0$, we would have from Lemma \ref{almostdecreasing}

    $$S(2t)-S(t) \le 8 \sqrt{N_0}.$$
    which is inconsistent with $S(t) \sim \frac{1}{2t^2 \zeta(\frac{1}{2})}$ as $t\rightarrow 0$, since  $\zeta(\frac{1}{2})<0$.

This does not rule out the possibility that $\sum_{n \le x} \frac{\lambda(n)}{\sqrt{n}}$ is negative for every large $x$, as is suggested in \cite{betweenpolyaandturan}.
\end{remark}
    
\bibliographystyle{plain}
	\bibliography{PD}{}

@article {kal2021,
    AUTHOR = {Kalmynin, A. B.},
     TITLE = {Quadratic characters with positive partial sums},
   JOURNAL = {Mathematika},
  FJOURNAL = {Mathematika. A Journal of Pure and Applied Mathematics},
    VOLUME = {69},
      YEAR = {2023},
    NUMBER = {1},
     PAGES = {90--99},
      ISSN = {0025-5793,2041-7942},
   MRCLASS = {11N05 (11K65 11L03)},
  MRNUMBER = {4516802},
MRREVIEWER = {Ben\ Joseph\ Green},
}

@incollection {BM89,
    AUTHOR = {Baker, R. C. and Montgomery, H. L.},
     TITLE = {Oscillations of quadratic {$L$}-functions},
 BOOKTITLE = {Analytic number theory ({A}llerton {P}ark, {IL}, 1989)},
    SERIES = {Progr. Math.},
    VOLUME = {85},
     PAGES = {23--40},
 PUBLISHER = {Birkh\"{a}user Boston, Boston, MA},
      YEAR = {1990},
   MRCLASS = {11M06},
  MRNUMBER = {1084171},
MRREVIEWER = {Matti Jutila},
}

@article{Jung2016,
  title   = {On the sparsity of positive-definite automorphic forms within a family},
  author  = {Jung, Junehyuk and Shin, Sug Woo},
  journal = {Journal d'Analyse Math{\'e}matique},
  volume  = {129},
  number  = {1},
  pages   = {105--138},
  year    = {2016},
  doi     = {10.1007/s11854-016-0017-9},
  url     = {https://doi.org}
}

@unpublished{Sarnak2011letter,
  author = {Sarnak, Peter},
  title = {Letter to {E}. {Bachmat} on positive definite {L}-functions},
  year = {2011},
  note = {Available at: \url{https://publications.ias.edu/sites/default/files/positive%20definite%20Lfunctions%20april.pdf}},
}

@article{aymone2024ii,
  author    = {Aymone, Marco},
  title     = {Sign changes of the partial sums of a random multiplicative function {II}},
  journal   = {Comptes Rendus. Math{\'e}matique},
  volume    = {362},
  pages     = {895--901},
  year      = {2024},
  doi       = {10.5802/crmath.615},
  url       = {https://doi.org}
}

@article{klurman2024signchangesshortcharacter,
      title={Sign changes of short character sums and real zeros of {F}ekete polynomials}, 
      author={Oleksiy Klurman and Youness Lamzouri and Marc Munsch},
      year={2024},
      note={Preprint available at arxiv.org/abs/2403.02195} 
}

@article{lamzouri2026realzeroslschid,
	 title={Real zeros of ${L}'(s, \chi_d)$}, 
      author={Youness Lamzouri and Kunjakanan Nath},
	year = {2026},
    note={Preprint available at arxiv.org/abs/2511.02774}
}

@article{kucheriaviy2025,
	title = {Positivity of partial sums of a random multiplicative function and corresponding problems for the Legendre symbol},
	author={Petr Kucheriaviy},
	year = {2025},
    note={Preprint available at arxiv.org/abs/2510.25691}
}

@article{Harpersup,
   title={Bounds on the suprema of Gaussian processes, and omega results for the sum of a random multiplicative function},
   volume={23},
   ISSN={1050-5164},
   url={http://dx.doi.org/10.1214/12-AAP847},
   DOI={10.1214/12-aap847},
   number={2},
   journal={The Annals of Applied Probability},
   publisher={Institute of Mathematical Statistics},
   author={Harper, Adam J.},
   year={2013},
   month=apr }

@article{FHK,
  title = {Freezing Transition, Characteristic Polynomials of Random Matrices, and the Riemann Zeta Function},
  author = {Fyodorov, Yan V. and Hiary, Ghaith A. and Keating, Jonathan P.},
  journal = {Physical Review Letters},
  volume = {108},
  issue = {17},
  pages = {170601},
  numpages = {5},
  year = {2012},
  month = {Apr},
  publisher = {American Physical Society},
  doi = {10.1103/PhysRevLett.108.170601},
  url = {https://link.aps.org}
}

@article{ASX,
      title={Nonnegative quadratic character sums}, 
      author={Rodrigo Angelo and K. Soundararajan and Max Wenqiang Xu},
      year={2026+},
     note={In preparation}, 
}

@article{fhk1,
      title={The Fyodorov-Hiary-Keating Conjecture. {I}}, 
      author={Louis-Pierre Arguin and Paul Bourgade and Maksym Radziwill},
      year={2020},
      note={Preprint available at arxiv.org/abs/2007.00988}, 
}

@article{fhk2,
      title={The Fyodorov-Hiary-Keating Conjecture. {II}}, 
      author={Louis-Pierre Arguin and Paul Bourgade and Maksym Radziwill},
      year={2023},
      note={Preprint available at arxiv.org/abs/2307.00982}, 
}

@article{betweenpolyaandturan, title={BETWEEN THE PROBLEMS OF {P}òLYA AND {T}URáN}, volume={93}, DOI={10.1017/S1446788712000201}, number={1–2}, journal={Journal of the Australian Mathematical Society}, author={Mossinghoff, Michael J. and Trudgian, Timothy S.}, year={2012}, pages={157–171}}

@article{lyapunovbound,
author = {Bentkus, V.},
title = {A Lyapunov-type Bound in Rd},
journal = {Theory of Probability \& Its Applications},
volume = {49},
number = {2},
pages = {311-323},
year = {2005},
doi = {10.1137/S0040585X97981123},

URL = { 
    
        https://doi.org/10.1137/S0040585X97981123
    
    

},
eprint = { 
    
        https://doi.org/10.1137/S0040585X97981123
    
    

}
}

\end{document}